%% file: article.tex
\documentclass{amsart}
\usepackage{graphicx}
\usepackage{tikz}

\usepackage[style=alphabetic, maxcitenames=50,maxnames=50,maxnames=5,minnames=3,maxbibnames=99]{biblatex}
\numberwithin{equation}{section}

\include{macros}
\title{Langlands duality for skein modules of 3-manifolds}
\author{David Jordan}
\date{June 2022}

\begin{document}
\begin{abstract}
I introduce new Langlands duality conjectures concerning skein modules of 3-manifolds, which we have made recently with David Ben-Zvi, Sam Gunningham, and Pavel Safronov.  I recount some historical motivation and some recent special cases where the conjecture is confirmed.  The proofs in these cases combine the representation theory of double affine Hecke algebras and a new 1-form symmetry structure on skein modules related to electric-magnetic duality.  This note is an expansion of my talk given at String Math 2022 in Warsaw, and is submitted to the String Math 2022 Proceedings publication.
\end{abstract}
\maketitle

\section{Introduction}
%Langlands duality conjectures are a cornerstone of contemporary research into arithmetic geometry, algebraic geometry, higher categories and quantum field theory.  Langlands duality patterns have been conjectured in \textbf{arithemetic} geometry and in the \textbf{de Rham} and \textbf{Betti} moduli spaces attached to algebraic curves.  Each of these conjectures has been \textbf{quantized}, and the quantum conjectures have found natural interpretation in \textbf{QFT}.

In his 1967 letter to Andr\'e Weil, Robert Langlands proposed a mysterious conjectural correspondence -- now known as Langlands reciprocity -- between what are now called \textbf{automorphic} representations of a simple algebraic group $G$ and \textbf{Galois} representations to its \emph{Langlands dual} group $^LG$.  Let us fix a number field $\mathbb{F}$ with ring of adeles $\mathbb{A}_\mathbb{F}$. Automorphic representations of $G$ are, loosely speaking, certain $G(\mathbb{A}_\mathbb{F})$-representations which may be realised inside $L^2(G(\mathbb{A}_\mathbb{F})/G(\mathbb{F}))$, and hence are closely related to automorphic forms.   The Langlands dual group $^LG$ is another simple algebraic group obtained from $G$ by interchanging root data.  Galois representations to $^LG$ are group homomorphisms from the absolute Galois group $\Gamma(\overline{\mathbb{F}}/\mathbb{F})$ to $^LG(\mathbb{C})$.  Among the many remarkable features of Langlands reciprocity already evident is that the objects it relates are of {\it a priori} very different nature. Langlands reciprocity and its many consequences and relatives -- including variants for function fields and local fields -- became collectively known as \textbf{arithmetic Langlands duality}.  %They have driven many advances in number theory for decades, but remain very much open.  

   \medskip

Decades after Langlands circulated his conjectures, Beilinson and Drinfeld discovered a new kind of Langlands duality -- also conjectural -- taking place between certain moduli spaces of bundles over a smooth projective curve $X$.  Their proposal, now known as the \textbf{de Rham geometric Langlands duality}, asserts an equivalence between different algebro-geometric categories living over two such moduli spaces.  On the automorphic side lies a certain category of $\D$-modules -- i.e. systems of polynomial differential equations -- on a moduli space of holomorphic $G$-bundles on $X$.  On the Galois side lies a category of coherent sheaves -- i.e., finite quotients of maps between vector bundles -- on a moduli space of flat $^LG$-bundles (i.e. $^LG$-bundles equipped with a flat connection 1-form) on $X$.  As in the arithmetic setting, the geometric Langlands duality relates objects of {\it a priori} very different nature, and bearing no elementary relationship.  

The thread connecting the arithmetic and geometric conjectures passes through a deep series of analogies -- known as Weil's Rosetta stone -- which relate the arithmetic of function fields to the geometry of complex curves.  One imagines the ring of adeles $\mathbb{A}_\mathbb{F}$ of a function field $\mathbb{F}$ to be a complex curve, with the primes as points, and hence regards $G(\mathbb{A}_{\mathbb{F}})$ as giving ``local coordinates" of a $G$-bundle in the formal neighborhood of each prime; on the geometric side this corresponds to the specification of a $G$-bundle by its Taylor expansion at each point of a curve.  One thereby situates the category of $\D$-modules on $\Bun_G$ on the geometric side opposite to a category of $\ell$-adic sheaves on the arithmetic side; these in turn are a kind of categorification of automorphic forms, via the ``sheaf-to-function" correspondence for function fields.  On the other side of Langlands duality, the absolute Galois group is well-understood to be an arithmetic analog of the fundamental group of the curve. \medskip 

%The drive to understand geometric Langlands conjectures has greatly furthered our understanding of stacks, $\infty$-categories, and $\mathcal{D}$-modules, and is a major motivator for the field of geometric representation theory.  Countless consequences and related dualities have been confirmed, and the geometric Langlands conjectures are generally now perceived to be within reach.\medskip

Beilinson and Drinfeld's conjectures were later deformed in the language of twisted $\D$-modules by Feigin, Frenkel and Gaitsgory, following a proposal of Stoyanovsky.  In their \textbf{quantum de Rham geometric Langlands duality} a new symmetry emerges which is not present in the arithmetic or geometric settings: the distinction between automorphic and Galois sides evaporates, and both moduli spaces which appear are those of holomorphic $G$-bundles (respectively, $^LG$-bundles), while the categories appearing are both categories of ($\Psi$- or ${^L\Psi}$-twisted, respectively) $\D$-modules.  Here the twisting depends on a possibly infinite complex parameter $\Psi\in\CP^1$, and ${^L\Psi}=-1/(n_G\Psi)$, where $n_G\in\{1,2,3\}$ is the lacing number of $G$.  The classical geometric Langlands conjecture is the special case $\Psi=\infty, {^L\Psi}=0$.

A physical manifestion of quantum geometric Langlands duality was discovered by Kapustin and Witten, as an instance of \textbf{$S$-duality}, also known as Montonen--Olive duality. This duality relates certain low-energy approximations called ``twists" of $\mathcal{N}=4$ super Yang-Mills quantum field theory in four dimensions.  $S$-duality is a sort of non-abelian Fourier transform, whose existence is predicted by $\mathcal{M}$-theory, and which can be understood as a generalisation of electric-magnetic duality in Maxwell theory.  In this telling $G$ and ${^LG}$ each appear as gauge groups of distinct theories, and the parameters $\Psi,{^L\Psi}$ each give $\mathbb{CP}^1$ charts on the space of twists.  $S$-duality asserts an equivalence of physical theories,
\[
\mathcal{Z}_{G,\Psi} \simeq \mathcal{Z}_{{^LG},{^L\Psi}},
\]
and Kapustin and Witten explained how to derive the categories appearing in both the classical and quantum geometric Langlands conjectures as categories of surface operators $\mathcal{Z}_{G,\Psi}(X), \mathcal{Z}_{{^LG},{^L\Psi}}(X)$.  This reformulation has been very influential in both mathematics and physics communities, and motivates much of the discussion which follows. \medskip

%The moduli spaces appearing in the de Rham geometric Langlands conjectures are notoriously technical in nature, owing to the highly infinite nature of de Rham connections, and of the group of gauge transformations:  the required moduli stacks are constructed as an infinite inductive limit of schemes, quotiented by an infinite inductive limit of groups.  In such settings the categories -- of $\mathcal{D}$-modules, of coherent sheaves -- appearing in the duality conjectures exhibit delicate technical behaviour even in their definition.  Indeed, much of the  foundational challenge around geometric Langlands conjectures is concerned with making precise sense of the categories involved. 

Ben-Zvi and Nadler proposed a further combinatorial incaranation of the geometric Langlands duality, now known as the \textbf{Betti geometric Langlands duality}.  The key idea is to apply the Riemann--Hilbert correspondence on both sides of de Rham conjectures.  On the Galois side, to a de Rham local system on $X$ is associated its monodromy homomorphism $\pi_1(X)\to {^LG}$.  On the automorphic side, to a $\cD$-module on $\Bun_G(X)$ is associated its sheaf of flat sections over $\Bun_G(X)$.  The Betti geometric Langlands duality conjectures assert an equivalence between a suitable category of coherent sheaves on the character variety -- the moduli space of homomorphisms $\pi_1(X)\to {^LG}$ -- and a category of sheaves of vector spaces on $\Bun_G(X)$.  In particular the latter is conjectured, like the former, to depend only on the topological surface $\Sigma$ underlying the algebraic curve $X$.\medskip

The Galois side of the Betti correspondence was subsequently $q$-deformed in the framework of quantum groups in our works with Brochier, Ben-Zvi, Snyder and Safronov, where it was called the quantum character theory TFT, and given the structure of a fully extended four-dimensional TFT (with divergent partition function, what is sometimes called a once-categorified 3-dimensional TQFT, or a (3+1)-TQFT).  The essential idea in these constructions was that the Betti Galois moduli spaces may be regarded as mapping stacks $\operatorname{Maps}(\Sigma,BG)$ into the classifying stack $BG$, with $\Coh(BG)=Rep(G)$.  The stack $BG$ has a canonical 2-shifted symplectic structure, whose deformation quantization $BG_q$ is given algebraically by the braided tensor category $\Rep_q(G)$.  Using the language of higher Morita theory and factorization homology, it was possible to make formal sense of a functor ``$\operatorname{Maps}(\Sigma,BG_q)$", and to establish that this extends to a 4-dimensional TQFT with divergent partition function.\medskip

\begin{center}\rule{.3\textwidth}{0.15ex}\end{center}\medskip

The purpose of this note is to share a recent conjecture we have made jointly with David Ben--Zvi, Sam Gunningham and Pavel Safronov, which envisions Langlands duality between the skein modules of closed oriented 3-manifolds.  The $G$-skein module $\Sk_{G,q}(M)$ of a closed, oriented 3-manifold is the formal linear span of $G$-labelled ribbon graphs, modulo local relations modelled on the ribbon braided tensor category $\Rep_q(G)$ of representations of the quantum group $U_q\mathfrak{g}$; skein categories are defined in a similar spirit (see Section \ref{sec:QGskeins} for precise definitions and some examples). In its simplest form, our conjecture states:

\begin{conjecture}\label{conj:main-conj}
Let $G$ be a semisimple algebraic group, and let ${^LG}$ denote its Langlands dual group.  Let $M$ be a closed, oriented 3-manifold, suppose that $\Psi \in \mathbb{C}^\times$ is transcendental, and let $q=e^{\mathrm{i}\Psi}$ and ${^Lq}=e^{\mathrm{i}{^L\Psi}}$.  Then we have a linear isomorphism,
\[
\Sk_{G,q}(M) \cong \Sk_{{^LG},{^Lq}}(M).
\]
\end{conjecture}

\noindent See Section \ref{sec:bigraded-refinement} for the statement of various refinements and strengthenings of the conjecture, and Section \ref{sec:evidence} for proofs in a few special families of examples.\medskip

\begin{remark}
It was conjectured by Edward Witten and proved in our work with Gunningham and Safronov that the skein modules are in fact finite-dimensional under the above assumptions, so that our conjecture can be rephrased as an equality of natural numbers:
\[
\dim\Sk_{G,q}(M) = \dim\Sk_{{^LG},{^Lq}}(M).
\]
Moreover our method of proof established that the dimension is independent of $q$ so long as it is taken to be transcendental (it is expected to suffice to require only that $q$ is not root of unity).  Hence the precise specification of parameters $q$ and ${^Lq}$ in the conjecture is somewhat superfluous, and is included only to give the reader a sense for the parallel with de Rham quantum geometric Langlands conjectures.

In particular, to confirm the conjecture as stated it is enough to compute the generic dimension of the skein module for $G$ and for ${^LG}$.  In fact, all the cases where we can confirm the conjecture are by independently computing both dimensions.  Because of the transcendental relationship between $q$ and ${^Lq}$, a canonical isomorphism of vector spaces (e.g. one intertwining the action of the diffeomorphism group on each side) is likely to depend analytically, rather than algebraically, on the parameter $q$.
\end{remark}

The following table summarises the many instances of Langlands duality recounted so far.  Note the symmetry shared between the quantum de Rham and quantum skein module formulations. \smallskip

{\small
\begin{center}\begin{tabular}{c|c|ccc}
    Shorthand & Object of study & Automorphic/A side && Galois/B side \\\hline
    Arithmetic & Number field $\mathrm{F}$ &  $\{V\subseteq L^2(G(\mathbb{A}_\mathrm{F})/G(\mathbb{F}))\}$ & $\leftrightarrow$ &  $\{\rho:\Gamma(\overline{\mathrm{F}}/\mathrm{F})\to {^LG}\}$ \\\hline
    Classical \\\hline
    de Rham & Complex curve $X$ & $\D(\Bun_{G}(X))$ & $\leftrightarrow$ & $\Coh(\Loc_{{^LG}}(X))$ \\
    Betti & Real surface $\Sigma$ & $\Shv(\Bun_{G}(\Sigma))$ & $\leftrightarrow$ & $\Coh(\Ch_{{^LG}}(\Sigma))$\\\hline
    Quantum \\\hline
    de Rham & Complex curve $X$ & $\D_\Psi(\Bun_G(X))$ & $\leftrightarrow$ & $\D_{{^L\Psi}}(\Bun_{{^LG}}(X))$\\
    Betti & Real surface $\Sigma$ & $\Shv_q(\Bun_G(\Sigma))$ & $\leftrightarrow$ & $\Coh_q(\Ch_G(\Sigma))$\\
    Skein & 3-manifold & $\Sk_{G,q}(M)$ & $\cong$ & $\Sk_{{^LG},{^Lq}}(M)$
\end{tabular}\end{center}}\smallskip

The most basic motivation for the conjecture is that the ${^LG}$-skein module of a 3-manifold $M$ specialises at $q=1$ to the algebra of functions on the ${^LG}$-character variety of $M$, and the skein category of a surface $\Sigma$, specialises at \mbox{$q=1$} to the ${^LG}$-character stack of a surface $\Sigma$.  Indeed, a ${^LG}$-labelled ribbon graph defines in a straightforward way a holomorphic function on the ${^LG}$-character variety, essentially by sending a ${^LG}$-local system $E$ to the trace of an associated bundle living along the support of the graph: these functions, and their deformations, are called Wilson loop observables.  It is then natural to hope -- though far from automatic -- that the one-parameter family of quantum deformations given by the skein module indeed coincides with the one-parameter family of twists in the Kapustin--Witten theories.  The rich structure of extended topological field theory enjoyed by skein modules -- they are expected to coincide with the value of the quantum character theory TFT on 3-manifolds, and hence to describe the Kapustin--Witten state space -- as well as their natural interpretation as $q$-deformed Wilson loop operators lends support for this hope.

Another source of motivation for our conjecture comes from the deep analogies of Mazur, Kapranov, and Reznikov, relating 3-manifold topology to algebraic number theory, and from the subsequent development of \emph{arithmetic quantum field theory}.   In this ``MKR dictionary" -- which builds on Weil's Rosetta stone -- one imagines a number field (more precisely its ring $\mathbb{O}$ of integers) to be a closed 3-manifold; one treats ideals in $\mathbb{O}$ as links, and prime ideals as knots.  Units in $\mathbb{O}$ are treated as embedded surfaces, field extensions are treated as branched covers, homology groups as ideal class groups, etc.  Building on this dictionary, the impetus of arithmetic field theory is to compute arithmetic invariants of number fields as if they were partition functions for a quantum field theory.  This idea was explored first by Kim, who formulated an arithmetic analog of 3D Chern-Simons TQFT as a mechanism for constructing $L$-functions, and taken up more recently by Ben-Zvi, Sakellaridis, and Venkatesh who have proposed to study various aspects of arithmetic Langlands duality using ideas from Kapustin and Witten's 4D TQFT. 

It is perhaps remarkable in hindsight given the MKR dictionary that our conjecture so significantly post-dates the geometric Langlands conjectures. Just as skein theory attaches a vector space to each 3-manifold, arithmetic Langlands duality attaches vector spaces -- the space of automorphic forms, and the space of algebraic functions on the arithmetic character variety, on each side of the duality -- to each number field.  The geometric Langlands dualities differ in two key respects: they involve complex curves ($\leftrightarrow$ function fields) as opposed to 3-manifolds ($\leftrightarrow$ number fields), and they replace vector space-valued invariants with categorical ones.  This of course makes Langlands duality for skein modules more elementary to formulate, and to falsify, than geometric Langlands duality, since it reduces to a statement about equality of integer dimensions rather than about equivalences of $\infty$-categories.  From a more physical perspective, it is interesting to note that Kapustin and Witten's expressed quantum geometric Langlands duality for \emph{2-dimensional manifolds} using a \emph{four-dimensional} QFT, but largely left aside the intervening question of mathematically describing the Hilbert spaces attached to \emph{3-manifolds}.

This sense of anachronism can be partly explained by posing the following two natural and (to our knowledge) unanswered questions concerning Langlands duality in the arithmetic and physical settings.  The answers to each question would fill in some more of the dots between our conjectures and the canon of Langlands duality which I have surveyed above.   On the side of arithmetic:

\begin{question}
What sense, if any, can be made of a \textbf{quantum arithmetic Langlands duality}?  
\end{question}

In particular, what role does the symplectic structure on character varieties, and its deformation quantization, play in number theory, and what is the arithmetic meaning of the deformation parameter $q$? On the physical side:

\begin{question}
What is the Hilbert space attached to a closed oriented 3-manifold by the Kapustin--Witten $A$-side twist at $\Psi=0$? 
\end{question}

The answer to this question would sit across Langlands duality from the ${^LG}$-character variety of $M$, and should have an $A$-side flavour involving symplectic/contact geometry, Fukaya categories, Floer theory, etc.  One may contemplate possible answers either by degenerating our understanding for generic $\Psi$ (as a skein module), or by extrapolating up in dimension from the case of surfaces, where the $A$-side twist at $\Psi=0$ is a category of $\mathcal{D}$-modules on $\Bun_G(X)$.

One among many complications which arise when contemplating this $A$-side twist at $\Psi=0$ for 3-manifolds is that while $S$-duality predicts an equivalence of theories $\mathcal{Z}_{G,\Psi} \simeq \mathcal{Z}_{{^LG},{^L\Psi}}$, it nevertheless \emph{interchanges} the Dirchlet and Nahm pole boundary conditions.  The skein module is, in some sense by definition, the orbit of the Dirichlet boundary condition by the Wilson loop operators of the theory; and for irrational non-zero $\Psi$ (equivalently, $q$ not a root of unity) we anticipate that these orbits will coincide.  We however do not expect such a coincidence at $\Psi=0,\infty$, and therefore one must formulate the Nahm pole boundary condition on 3-manifolds in mathematical terms.

\subsection*{References} Many of the conjectures recalled in this introduction were circulated first as private correspondence, were never published, or were published only much later as folklore.  For this reason, I have avoided in-line citations throughout the introduction, even in those cases which appeared promptly in print, such as Kapustin and Witten's paper \cite{KapustinWitten} on electric-magnetic duality and Ben-Zvi and Nadler's paper \cite{BZN18} introducing Betti geometric Langlands conjectures.  For a more thorough recounting of Langlands duality in its arithmetic, geometric, and physical manifestions, and for an introduction to arithmetic topology, I recommend:

\begin{itemize}
    \item \emph{The genesis of the Langlands program} \cite{Langlandshistory}, an LMS Lecture series which covers Robert Langlands's biography, his mathematics, and his letter to Andre Weil,
    \item \emph{Between electric-magnetic duality and the Langlands program} \cite{BZlec}, an unpublished series of lecture notes, on the relations between geometric Langlands, arithmetic Langlands, and electric-magnetic duality.
    \item \emph{Knots and primes} \cite{knotsandprimes}, a Springer textbook introducing the foundations of arithemetic topology.
\end{itemize}

\subsection*{Acknowledgements}  I thank my collaborators David Ben-Zvi, Sam Gunningham, Pavel Safronov, Monica Vazirani, and Haiping Yang for our many discussions surveyed here, and for their encouragement to share our conjectures and work in progress in this forum.  I also thank Dan Freed, Du Pei and Edward Witten for several enlightening discussions and correspondences regarding electric-magnetic duality, and Minhyong Kim for inspiring discussions about arithmetic quantum field theory.  Finally, I thank the organisers of String Math 2022 for a great conference, and for the opportunity to prepare these notes.

\section{Quantum groups and their skein theory}
\label{sec:QGskeins}
%Let us now recall the formal definition of skein modules, their basic properties, and also their relation to character varieties, when specialised to $q=1$.

The construction of skein modules and skein categories depends on the specification of a cp-rigid ribbon braided tensor category\footnote{For experts:  equivalently, we may consider instead the subcategory of compact projective objects, which will be a (not-necessarily-unital!) ribbon braided monoidal tensor category in which all objects are dualizable.  We note that unless $\cA$ is semisimple, the compact-projective objects do not form an abelian category, but instead an additive idempotent complete category.}  $\cA$.  This means that $\mathcal{A}$ is a locally presentable abelian category with enough compact projective objects, equipped with a bilinear monoidal functor $\otimes: \mathcal{A} \times \mathcal{A} \to \mathcal{A}$, a braiding natural isomorphism $\sigma$, and a ribbon element $\nu$, in which all compact projective objects are dualizable.

The idea will be that we can encode all the algebraic structure of a cp-rigid ribbon braided tensor category in a diagrammatic form, using the \emph{graphical calculus} of Reshetikhin and Turaev.  These diagrams and all relations between them take place in a decorated cylindrical ball $D^2\times I$, but we can then allow our diagrams to permeate an entire oriented 3-manifold, imposing only locally the relations coming from the graphical calculus in some ball.  The resulting vector space is called the skein module.  We may additionally allow these diagrams to end at boundary surface components of our 3-manifold; this leads to the notion of a skein category, and of the induced functor between skein categories.

\subsection{Typical examples of ribbon braided tensor categories}

Our main examples of interest will be the categories $\Rep_q(G)$ attached\footnote{More precisely the specification of ribbon structure depends on choosing a $d$th root of $q$, where $d$ is the determinant of the Cartan matrix, but we will suppress this point from the notation.} to an algebraic group $G$ at transcendental values of the parameter $q$.  This means we consider the category of locally finite-dimensional representations of the quantum group $U_q(\mathfrak{g})$ attached to the Lie algebra $\mathfrak{g}$ of $G$, and within it we consider a tensor subcategory consisting of those $U_q(\mathfrak{g})$-representations, all of whose weights lie in the weight lattice of $G$.  We recall that the weight lattice of $G$ lies between those of the adjoint and the simply connected group of $\mathfrak{g}$, and hence is determined by specifying $Z(G)$, the quotient of the weight lattice of $G$ by the root lattice; or equivalently by $\pi_1(G)$, the quotient of the coweight lattice by the coroot lattice of $G$.

\begin{example}\label{ex:SL2PGL2}
We will treat  the examples $\Rep_q(\SL_2)$ and $\Rep_q(\PGL_2)$ in detail throughout these notes. We recall that $\Rep_q(\SL_2)$ is a semisimple category with simples $\mathrm{V}(0), \mathrm{V}(1), \mathrm{V}(2), \ldots$ with each $\mathrm{V}(k)$ having dimension $k+1$ and highest weight $k$.  The category $\Rep_q(\PGL_2)$ is the full subcategory of $\Rep_q(\SL_2)$ consisting of direct sums of the irreducible representations $\mathrm{V}(0), \mathrm{V}(2), \mathrm{V}(4), \ldots$.

Let us recall some further standard algebraic facts about these two categories.  First, in the case of $\Repq(\SL_2)$, we recall that the defining two-dimensional representation $\mathrm{V}(1)$ is a \emph{tensor generator}: every representation $\mathrm{V}(k)$ appears as a direct summand in the kth tensor power $\mathrm{V}(1)^{\otimes k}$.  Secondly, we have that the $\Hom$ spaces $\Hom(\mathrm{V}(1)^{\otimes k}, \mathrm{V}(1)^{\otimes \ell})$ are freely generated by the evaluation and coevalution maps depicted in the first term of the RHS below.  In particular, we can express the braiding in terms of the evaluation, coevaluation and identity maps:\smallskip

\begin{center}\includegraphics[height=4cm]{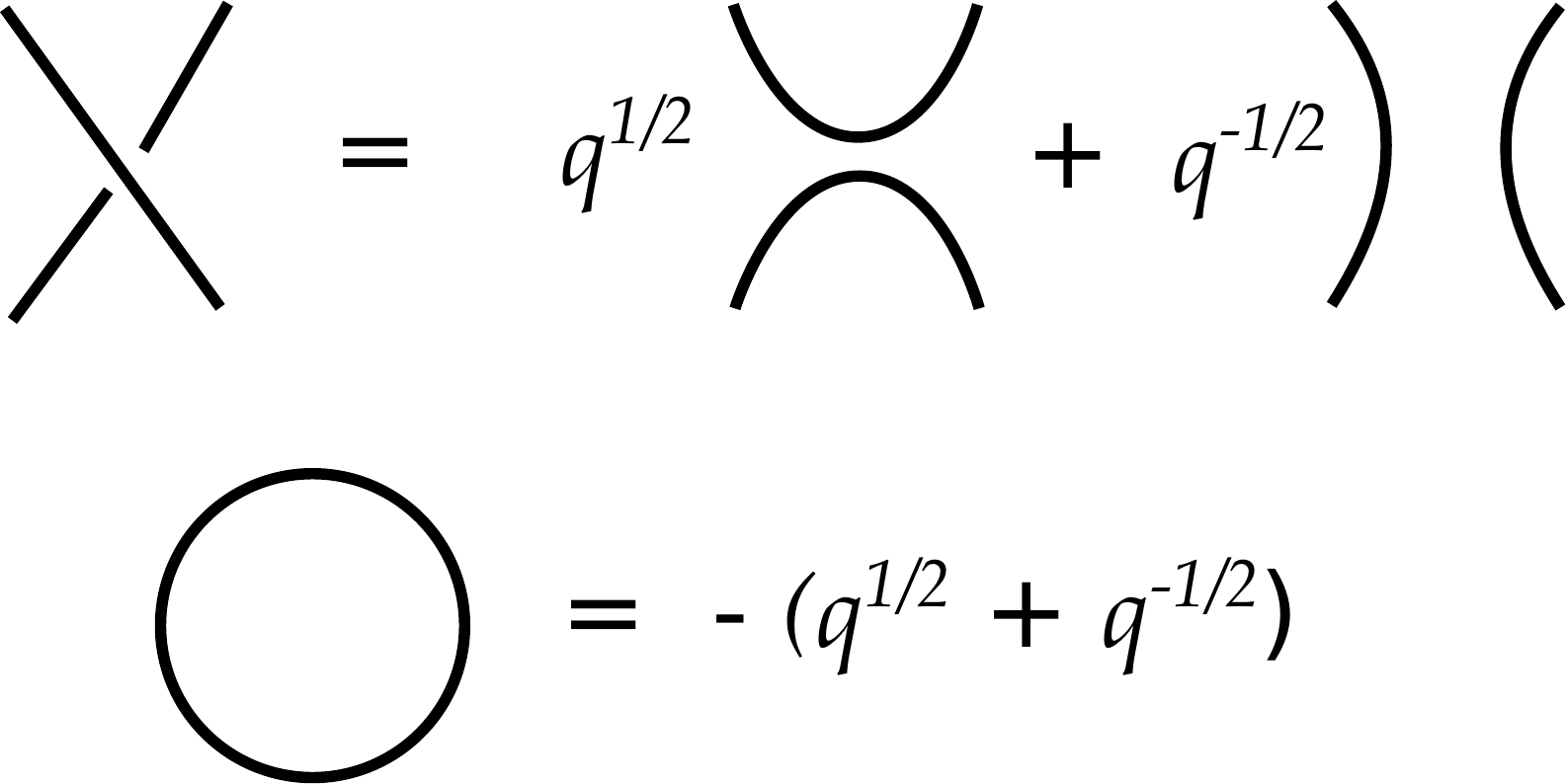}\end{center}\smallskip

We have omitted the source and target $\mathrm{V}(1)$ from the generators.  Let us note that the evaluation and coevaluation ought usually to involve $\mathrm{V}(1)^*$, but that we have identified $\mathrm{V}(1)\cong \mathrm{V}(1)^*$.  There are two natural candidates for such isomorphisms, and it is customary to favor the skew-symmetric pairing, hence the negative sign in the formula for the value of the circle, i.e. the quantum dimension of $\mathrm{V}(1)$.
\end{example}

\begin{example}\label{ex:PGL2skein}
For $\Repq(\PGL_2)$, the situation is similar, but somewhat more complicated.  It is again easy to show that $\Repq(\PGL_2)$ is tensor generated by the generator $\mathrm{V}(2)$.  Recalling the isomorphism $\mathrm{\mathrm{V}}(2)\otimes\mathrm{V}(2) \cong \mathrm{V}(0) \oplus \mathrm{V}(2) \oplus \mathrm{V}(4)$, we have now a three-dimensional space $\operatorname{End}(\mathrm{V}(2)\otimes \mathrm{V}(2))$, spanned by the three diagrams on the RHS below.  In particular, we may express the braiding:\smallskip

\includegraphics[height=2cm]{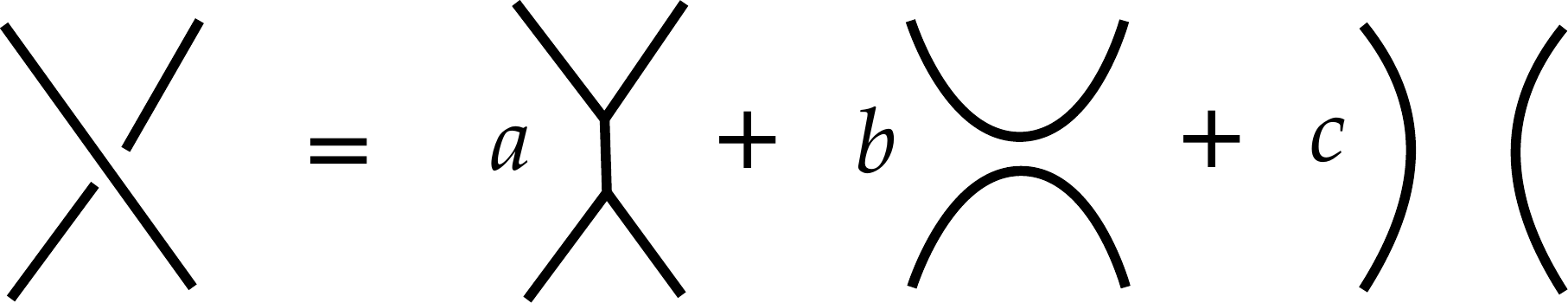}\smallskip

\noindent for suitable constants $a, b, c \in \C(q)$ which can be determined by direct computation in the given basis. Here again we have omitted the label $\mathrm{V}(2)$ from the generators.
\end{example}

\begin{example}
For the family of groups $G=\SL_N$, we also have that $\Repq(\SL_N)$ is tensor generated by the defining representation $\mathrm{V}(\omega_1)$ associated to the first fundamental weight of $\SL_N$.  It is possible to give a similar explicit diagrammatic description for $\Hom$ spaces of tensor powers of $\mathrm{V}(\omega_1)$, using either MOY calculus \cite{MOY} or the spiders/webs descriptions of \cite{Kuperberg,CKM}, to which we refer for details.
\end{example} 

\subsection{Graphical calculus, skein modules and skein categories}
The key fact we will use about ribbon braided tensor categories is that they admit a \emph{graphical calculus}.  We have already seen some examples of graphical calculus above; let us recall what it means more precisely. 

\begin{figure}
\includegraphics[height=4cm]{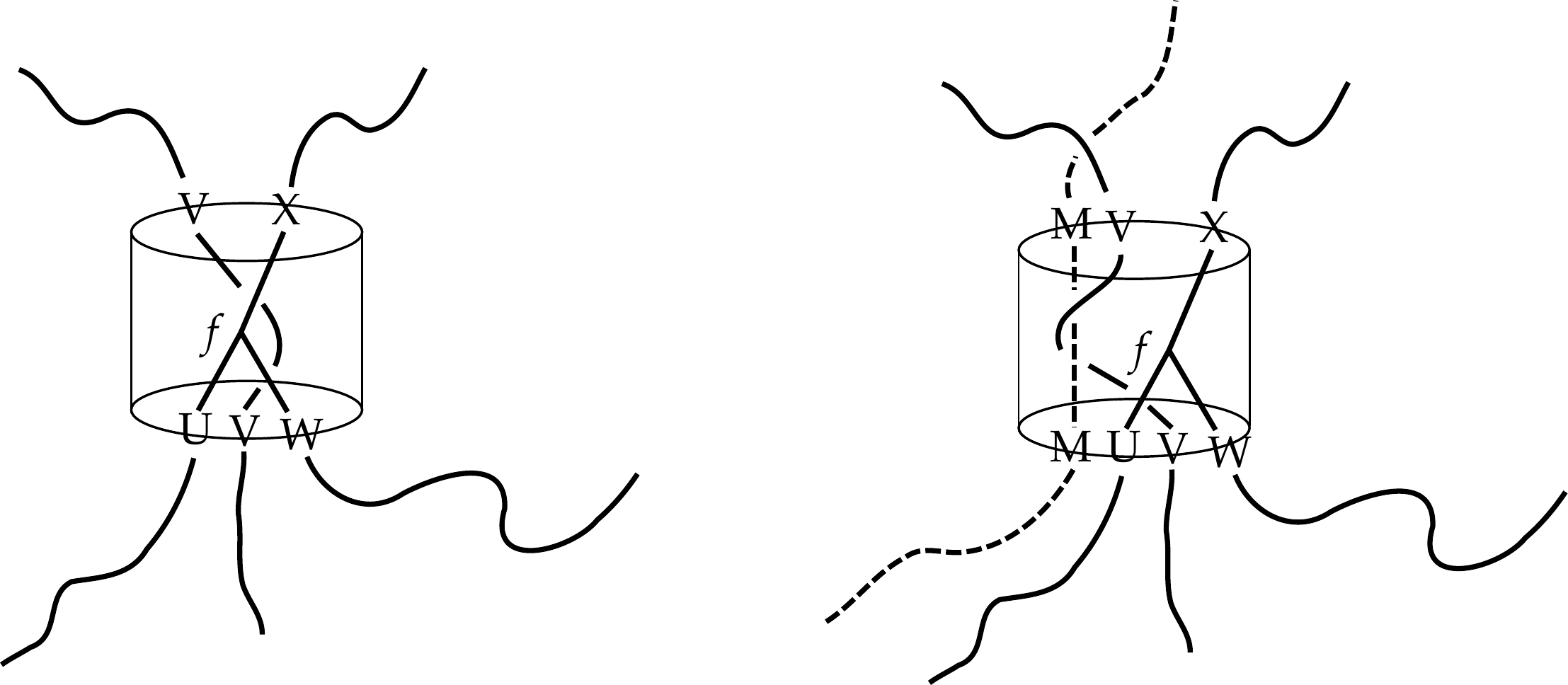}
\caption{At left: the skein module is spanned by $\cA$-labelled ribbon diagrams; relations; the skein relations state that the diagram within the cylinder can be replaced by any other diagram whose RT-evaluation gives the same element of $\Hom(U\otimes V \otimes W,V\otimes X)$.  At right: a defect is introduced along the dotted line, and a braided module category is inserted along it (see Section \ref{sec:gradingstwistings}).  The skein module now includes labels and relations from the braided tensor category away from the defect, and from its braided module category along the defect.}\label{fig:ribbon}
\end{figure}

\begin{definition}[Sketch.  See Figure \ref{fig:ribbon} and {\cite[Section 4.2]{Cooke}}]
	Let $\Sigma$ be an oriented surface.
	\begin{itemize}
		\item An \defterm{$\cA$-labeling of $\Sigma$} is the data, $X$, of an oriented embedding of finitely many disjoint disks $x_1, \dots, x_n\colon \DD\rightarrow \Sigma$ labeled by objects $V_1, \dots, V_n$ of $\cA$.  We denote by $\vec{x_i}$ the $x$ axis sitting inside each disk $x_i$, and denote $\vec{X}=\cup_i\vec{x_i}$.
		\item A ribbon graph has ``ribbons" connecting ``coupons".  A coupon is simply a zero-simplex in $\Sigma\times I$, equipped with the choice of a normal vector.  A ribbon is an embedded 1-simplex in $\Sigma\times I$, whose endpoints each lie either in the set $X$, or at some coupon, and which is equipped with a normal vector along it.  We furthermore require that the tangent vector of each embedded simplex points directly up when ending at $X\times \{0,1\}$ or else in the direction of the coupon's normal vector when ending at a coupon. Finally, we require a cyclic ordering of the incoming and outgoing ribbons at each coupon.
		\item An \defterm{$\cA$-coloring} of a ribbon graph is a labelling of each ribbon by an object of $\cA$, and of each coupon by a morphism from the (ordered) tensor product of incoming edges to the (ordered) tensor product of outgoing edges.
		\item We say that an $\cA$-colored ribbon graph $\Gamma$ is \defterm{compatible with an $\cA$-labeling} if $\partial \Gamma=\vec{X}$, and denote by $\mathrm{Rib}_\cA(M,X)$ the vector space with basis the $\cA$-colored ribbon graphs on $M$ compatible with $X$.
		\end{itemize}
	\end{definition}

\begin{remark} The graphical calculus for ribbon tensor categories has appeared in the literature with many differing but essentially equivalent formulations.  The formulation we present here is adapted from our presentation in \cite{GJS}, however the customary ``ribbons" and ``coupons" formulation has been replaced with a topologically equivalent formulation which will be more compatible with simplicial homology arising later in the paper.\end{remark}

Let us now describe the defining relations of skein modules.  Consider the 3-ball written as the cylinder $\DD\times I$ over the 2-disk $\DD$, and fix a labeling $X\cup Y$ with disks $X=(x_1,V_1),\ldots,(x_n,V_n)$ embedded in $\DD\times\{0\}$ and $Y=\{(y_1,W_1),\ldots (y_m,W_m)\}\times \{1\}$.  Then we have a well-defined surjection, 
	\[\mathrm{Rib}_\cA(\DD\times I,X\cup Y) \to \Hom_\cA(V_1\otimes \cdots \otimes V_n,W_1\otimes \cdots \otimes W_m),\]
	called the Reshetikhin--Turaev evaluation map, see \cite{TuraevBook}.  In essence this map takes a ribbon graph in a cylinder, and a generic projection, and interprets each crossing as a braiding, and each coupon as a morphism between each tensor product, and computes the corresponding composition: that this is well-defined independent of crossings is perhaps the most fundamental algebraic property of ribbon categories.  We will call the kernel of this map the \defterm{skein relations} between $X$ and $Y$.
	
	\begin{definition}
		Let $M$ be an oriented 3-manifold equipped with a decomposition of its boundary $\partial M\cong \rev{\Sigma}_{in}\coprod \Sigma_{out}$, and $\cA$-labelings $X_{in}$ of $\Sigma_{in}$ and $X_{out}$ of $\Sigma_{out}$. \begin{itemize} \item The \defterm{relative $\cA$-skein module} $\Sk_\cA(M, X_{in}, X_{out})$ is the vector space spanned by isotopy classes of $\cA$-colored ribbon graphs in $M$ compatible with $X_{in}\cup X_{out}$, taken modulo isotopy and the skein relations between $X_{in}$ and $X_{out}$ determined by any oriented ball $\DD\times I\subset M$\footnote{Here we assume without loss of generality that $\DD\times\{0\}\subset \Sigma\times \{0\}$ and $\DD\times\{1\}\subset \Sigma\times \{1\}$.}.
			\item When $\partial M=\emptyset$ (hence $\partial\Gamma=\emptyset$), we call this the \defterm{$\cA$-skein module}, and denote it by $\Sk_\cA(M)$.
		\end{itemize}
	\end{definition}
	
	\begin{definition}
		Let $\Sigma$ be an oriented surface. The \defterm{skein category} $\SkCat_\cA(\Sigma)$ of $\Sigma$ has:
		\begin{itemize}
			\item As its objects, $\cA$-labelings of $\Sigma$.
			
			\item As the 1-morphisms from $X$ to $Y$ the relative $\cA$-skein module of $(\Sigma\times [0, 1],X,Y)$.
		\end{itemize}
		The composition of 1-morphisms is by concatenation of relative skeins in the evident sense.
	\end{definition}
	
\begin{example}
In the cases $\cA=\Repq(\SL_2)$ and $\cA=\Repq(\PGL_2)$ the rather intricate definition of the skein module involving labelled ribbon graphs simplifies considerably, owing to the simplifications highlighted in Example \ref{ex:SL2PGL2}.  In both cases, one can dispense with labels and orientations on strands: for labels, one takes $\mathrm{V}(1)$ or $\mathrm{V}(2)$, respectively as the label always (and thus omits it), and one further applies the self-duality of $\mathrm{V}(1)$, respectively $\mathrm{V}(2)$ to eliminate orientation data.

For coupons, in the $\Repq(\SL_2)$ case, Schur Weyl duality implies that all morphisms can be expressed in terms of the braidings, so no coupons are necessary at all.  For $\Repq(\PGL_2)$, we require only the trivalent coupon from Example \eqref{ex:PGL2skein}.
\end{example}

\subsection{Relations to deformation quantization}
The starting point for studying skein modules using geometric representation theory comes from the following description of skein modules and skein categories at $q=1$.

\begin{definition}
The $G$-character stack $\underline{\Ch}_G(X)$ of a topological space $X$ is the moduli stack of flat $G$-bundles on $X$, equivalently, it is the quotient stack obtained by quotienting the \defterm{$G$-representation variety} $\Hom(\pi_1(X),G)$, by the conjugation action of $G$.  The \defterm{$G$-character variety} $\Ch_G(X)$ is the GIT quotient, i.e. the set of closed $G$-orbits on the representation variety.
\end{definition}

\begin{theorem}
The idempotent completion (respectively, free cocompletion) of the skein category $\SkCat_{G,q=1}(\Sigma)$ is naturally equivalent to the category of vector bundles (respectively quasi-coherent sheaves) on the $G$-character stack $\underline{\Ch}_G(\Sigma))$.
Moreover, the deformation to the skein category is flat in the parameter $q$.
\end{theorem}

Here, flatness in the parameter $q$ means that $\Hom$ spaces defined in $\SkCat_G(\Sigma)$ are actually free $\C[q,q^{-1}]$-modules; informally, one says then that the vector space supporting the skein module is ``constant" in $q$ and that only the algebra structure changes as a function of $q$: in such cases, we obtain a Poisson bracket on the $q=1$ limit, and we say that the skein algebra/skein category is a deformation quantization of the character variety/character stack.  The resulting Poisson bracket quantizes the Atiyah--Bott/Goldman/Fock-Rosly Poisson bracket (a single Poisson bracket on the character variety independently constructed by each grouping of authors in three different ways).  One has a similar interpretation for the $q=1$ limit of $\Sk_G(M)$ for any 3-manifold.  The following is proved for $G=\SL_2$ in \cite{BullockRings,PrzytyckiSikora}, but holds for any group $G$.

\begin{theorem}
The skein module $\Sk_{G,q=1}(M)$ of an oriented 3-manifold $M$ is naturally isomorphic to the algebra $\mathcal{O}(\Ch_G(M))$ of algebraic functions on the $G$-character variety of $M$.
\end{theorem}

However, when we move the quantum parameter $q$ away from one, skein modules behave very differently than skein categories in elementary terms.  While skein categories of surfaces (more precisely the relative skein algebras which appear as their $\Hom$ spaces) depend \emph{flatly} on the quantum parameter $q$, quite the opposite is true for skein modules of closed 3-manifolds.  We have noted that at $q=1$ the skein module of a closed 3-manifold may be identified with the algebra of algebraic functions on the $G$-character variety of $M$ -- this is in particular typically an infinite-dimensional vector space.  By contrast, we have

\begin{theorem}[\cite{GJS}]
Let $q\in \mathbb{C}^\times$ be generic, and let $M$ be a closed, oriented 3-manifold.  Then the $G$-skein module of $M$ is finite-dimensional.
\end{theorem}

This statement first appeared as a conjecture of Edward Witten in 2014, and came largely as a surprise to the skein theory community:  because of the many known results asserting that skein algebras and skein categories of surfaces deform flatly, it was generally expected that skein modules of 3-manifolds should as well.  With the benefit of hindsight and modern approaches to deformation quantization -- specifically the framework of shifted symplectic structures \cite{PTVV} -- the finiteness asserted in the theorem is not at all surprising.

The character stacks of surfaces, which are deformed by skein categories, are symplectic (what is now called ``0-shifted symplectic"), where as the character stacks of 3-manifolds have a very different structure: they are what is now called $(-1)$-shifted symplectic.  One way to understand $(-1)$-shifted symplectic stacks is by a shifted analog of Darboux coordinates:  whereas Darboux coordinates equip a 0-shifted symplectic stack with a local symplectomorphism to the unique symplectic vector space in that dimension, Darboux coordinates for $(-1)$-shifted stack identify the neighborhood of any point with the intersection of two Lagrangian subspaces of the symplectic vector space.

The finite-dimensionality of skein modules can be understood in this context from the fact that the notion of deformation quantization for 0-shifted and $(-1)$-shifted symplectic structures is very different in nature.  Whereas the familiar notion for $0$-shifted structures involves a non-commutative algebra, with  a flatness criterion on the underlying vector space of the quantization, the quantization of a $(-1)$-shifted symplectic structure defines a perverse sheaf on the classical variety, whose global sections describes the deformation.  There are (at least) two general avenues to construct this perverse sheaf: an algebraic construction of deformation quantizations and their tensor products, following Kashiwara and Schapira \cite{KashiwaraSchapiraDQ}, or a geometric construction following \cite{BBBBJ}.

In the algebraic construction, we suppose that the $(-1)$-shifted structure really does arise as the intersection of two Lagrangian subvarieties $L_1, L_2$, and that the symplectic variety $X$ has a deformation quantization into some non-commutative algebra $A$ (or sheaf of algebras), and that the Lagrangians each quantize to holonomic (left, and right, respectively) modules $M_1$, $M_2$  over $A$.  In this case, Kashiwara and Schapira show that the non-commutative intersection -- i.e. the relative tensor product $M_1\otimes_A M_2$ -- defines a perverse sheaf over the classical intersection $L_1\cap L_2$.  A careful analysis of the sections of this perverse sheaf in the setting of skein modules yielded the finiteness result of \cite{GJS}.

A related geometric construction of Abouzaid and Manolescu \cite{AbouzaidManolescu} of complexified instanton Floer homology also hinges on the construction of the character stack of a 3-manifold, as an intersection of two standard Lagrangians in the character stack of a surface of some genus $g$.  The Brav--Bussi--Ben-Bassat--Joyce quantization of such a $(-1)$-shifted symplectic stack is controlled by a perverse sheaf of \emph{vanishing cycles}, and in good situations this will have, indeed, a finite dimensional space of global sections.

\begin{remark}
Work in progress by Gunningham and Safronov gives an identification of the Kashiwara--Schapira sheaf and the Brav--Bussi--Ben-Bassat--Joyce sheaf in great generality, including in our setting of skein modules, thus unifying the two approaches.
\end{remark}

We might be tempted to understand the finiteness obtained in these two constructions as follows:  generically, two Lagrangian submanifolds should intersect transversely, but of course in practice they may not.  However, when we quantize, even if the classical Lagrangians fail to intersect transversely -- giving a positive dimensional character variety, and hence an infinite-dimensional algebra of functions -- at the quantum level, and at generic parameters, they exhibit the generic behavior.

Both the Kashiwara--Schapira and the Abouzaid--Manolescu quantizations are by construction not only vector spaces but chain complexes; we have so far only constructed on their cohomology in degree zero.  This suggests that the skein module also has a natural derived enhancement, with finite-dimensional cohomologies in each degree, and satisfying a Langlands duality akin to \ref{conj:main-conj}.  We will not discuss this derived enhancement in the present note.

\subsection{Skein topological field theory}
Given a 3-manifold $M$ with oriented boundary $\partial M = \Sigma_+ \sqcup \Sigma_-$, we obtain a functor,
\[\Sk_\cA(M): \SkCat(\Sigma_-)^{op}\otimes\SkCat(\Sigma_+)\to \Vect,\]
called the skein bimodule of $M$, which assigns to a pair of skein objects the relative skeins in $M$ compatible with the given marking at the boundary.

\begin{theorem}[Walker, unpublished]
The assignments of surfaces and 3-manifolds to their skein category and skein (bi)modules, respectively, defines a functor from the category $\operatorname{Cob}^3$ of surfaces and 3-cobordisms, to the category of categories and categorical bimodules.
\end{theorem}

%In the typical examples of interest, where $\cA=\Rep_q(G)$ for some generic complex parameter $q$, we will abbreviate notation and write $\Sk_G(M)$, $\SkCat_G(\Sigma)$, respectively for the skein module and skein category.

It was shown in \cite{Cooke} that skein categories satisfy a universal property called excision, which implies that the skein category coincides with the categorical factorization homology of Ayala--Francis \cite{AyalaFrancis}, as was studied in this context in \cite{BZBJ1,BZBJ2}.  The excision property for skein modules should be regarded as expressing the existence of an extension of skein theory all the way down to the point:  the $(3,2)$-TFT axioms allow us to cut 3-manifolds along surfaces, and excision allows us to cut surfaces up along embedded 1-manifolds, and then those 1-manifolds up along 0-manifolds.

However, while it is known that skein modules define a $(3,2)$-TFT after Walker, and that excision for skein categories implies a $(2,1,0)$-TFT via factorization homology, I am not aware of a rigorous construction of the desired $(3,2,1,0)$-TFT in the framework of fully extended topological field theories.  In \cite{BJS}, we have constructed a $(3,2,1,0)$-TFT from the same input, of a cp-rigid braided tensor category, however it remains to identify the value of this theory on 3-manifolds with the skein module, as one expects.  We will nevertheless largely pass notions between the skein-theoretic and fully extended formulations of the theory in Sections \ref{sec:gradingstwistings}.

\section{Gradings and twists of skein modules}\label{sec:gradingstwistings}
It is an elementary observation that the ribbon braided tensor category $\Rep_q(G)$ is naturally graded (as an abelian monoidal category) by the Pontryagin dual $Z(G)^\vee$ of the centre $Z(G)$ of the group $G$: the central elements of $G$ define grouplike central elements of the quantum group $U_q(\mathfrak{g})$, and by Schur's lemma, any simple representation uniquely determines a character of the centre.  The group-like property implies that the representations tensor compatibly with the group structure.

Because the skein relations all take place locally in a ball, the central grading on $\Rep_q(G)$ extends naturally to a grading by $H_1(M,Z(G)^\vee)$ on the $G$-skein module: given a $G$-labelled ribbon graph, we obtain a simplicial 1-cycle by replacing the label along each 1-simplex with its grading by $Z(G)^\vee$.  We note that we always obtain a cycle in this way, since at coupons, there can only be a nonzero morphism if the sum of degrees is zero. If the manifold has boundary, and the skein is defined relative to the boundary, then the 1-cycle will give a relative $Z(G)^\vee$ 1-cycle, by similar reasoning.  This assignment gives the grading of the $G$-skein module by $H_1(M,Z(G)^\vee)$, and this grading will play an important role in what follows.

We will see that skein modules carry another structure: they can be twisted by braided module categories inserted along codimension-two defects.  The grading and twisting constructions are dual in a precise sense, which we will soon see.  For now let us outline the construction of twisting of skein modules:  the complete details are presented in the forthcoming \cite{GJS2}, so we will be somewhat brief.

\begin{definition}\label{def:brmodcatdef}
A braided module category $\mathcal{M}$ for a braided tensor category $\cA$ is a module category $\cM$ for $\cA$, together with a `monodromy'' natural automorphism,
\[
\alpha: (X,m)\mapsto (\alpha_{X,m}: X\otimes m \to X\otimes m), 
\]
of the action functor of $\cA$ on $\cM$, which satisfies:
\[
\alpha_{X\ot Y,m} = (\sigma_{Y,X}\sigma_{X,Y}\ot \id_M)\circ \alpha_{X,Y\ot M}\circ (\id_X \otimes \alpha_{Y,M}).
\]
\end{definition}

\begin{remark}To motivate this definition, we should imagine that the braided tensor structure on $\cA$ is encoded by the $E_2$-operad of little disks (see Figure \ref{fig:brmodcat}), and that $\mathcal{M}$ sits at a single point (a co-dimension two defect) in the plane.  The module action comes from bringing $\cA$-disks to $\cM$, and the natural isomorphism $\alpha$ is indeed the monodromy obtained by sliding those disks once counterclockwise around $\cM$ before colliding.\end{remark}

\begin{definition}
The $E_2$-monoidal bicategory $\BrMod(\cA)$ has as objects braided module categories, as its 1-morphisms braided module functors (i.e. module functors $(F,J)$ satisfying an evident compatibility with $\alpha$), and as 2-morphisms the braided module natural transformations.  The $E_2$-monoidal structure is given by sending the pair $\cM_1$, $\cM_2$ to the relative tensor product $\cM_1\otimes_\cA \cM_2$, with the resulting holonomy natural automorphism given as the composition of holonomies.
\end{definition}

We will say that the braided module category $\cM$ is dualizable if it is fully dualizable as an object of $\BrMod(\cA)$.

Given a fully dualizable braided module category, and a knot $K$ in a 3-manifold $M$, we may define the skein module, $\Sk_{\cM}(M,K)$ as the formal span of all $(\cA,\cM)$-labelled ribbon graphs, modulo $(\cA,\cM)$-skein relations.  Here, an $(\cA,\cM)$-labelled ribbon graph is a ribbon graph in $M$, equipped with a labelling by $\cA$ in the open set $M\backslash K$, and with a labelling by $\cM$ along $K$.  The $(\cA,\cM)$-skein relations are given by imposing within any open cylinder contained in $M\backslash K$ the Reshetikhin--Turaev evaluation relations of $\cA$, and within any open cylinder intersecting $K$ along $\{0\}\times I \subset D^2\times I$ the Reshetikhin--Turaev evaluation relations of $\cM$.

\subsection{Higher form symmetry actions of abelian groups}

The most important source of dualizable module categories for the present discussion are the \emph{invertible} module categories, those for which the evalution and coevaluation are equivalences of braided module categories, identifying the dual braided module category as an inverse.  These form an abelian group under tensor product, and we can use this additional structure to define twists not only along knots but along simplicial cycles in $M$.  To motivate the construction, we recall now the framework of higher form symmetries in quantum field theory.

Recall that an ordinary action of an abelian group $A$ on a $(d+1)$-dimensional QFT induces an action of $H_{d}(M,A)$ on the Hilbert space $Z(M)$ attached to any $d$-manifold $M$ -- in particular if $M$ is closed and connected we have an action of the group $A\cong H_d(M,A)$ itself, on the vector space $Z(M)$.

In purely local/field-theoretic terms, the $A$-action arises because we redefine the field equations in the neighborhood of some $d$-cycle, requiring that fields transform according by the $A$-action as they pass through the cycle: in particular, this introduces a singularity to the field and so a field in the twisted sense is not a field in the ordinary sense.  The action of $H_d(M,A)$ on $Z(M)$ which we highlighted arises by situating a $d$-cycle at the halfway point of trivial mapping cylinders $M\times I$, thereby modifying the notion of correlation functions.

In the literature on extended defects, such an action is called a 0-form symmetry, the ``0" expressing that it is not yet an extended symmetry.  In this context, we may alternatively contemplate the structure of $k$-form symmetry action of a group $A$ on a $(d+1)$-dimensional QFT: in mathematical/functorial terms this means we consider an action of the $(k+1)$-group $B^k A$ on the theory (here $B^k$ denotes the $k$-fold delooping, a $k+1$-groupoid with unique object, and unique $n$=morphism for $n<k$ and with $A$ as its group of $k+1$-morphisms).  From such an extended action we obtain (among other things) an action of $H_{d-k}(M,A)$ on the Hilbert space $Z(M)$ attached to each $d$-manifold $M$.  Again in physical terms, the data of a $k$-form symmetry action of the group $A$ allows for modifications of the field equations, now along higher co-dimension submanifolds of $M$.

\subsection{Twisting and grading together: the Picard groupoid}

The homology groups $H_{d-k}(M,A)$ appearing above are \emph{global} invariants of $M$ -- that is, they depend on the global topology of $M$, and they allow us to twist Hilbert spaces globally along some global cycle.  It will be very useful in practice to give a more local enhancement of a $k$-form symmetry action, by introducing \emph{twists} $Z(M)[\gamma]$ of the Hilbert space $Z(M)$ indexed by elements $\gamma\in H_{d-k-1}(M,A)$, and isomorphisms between homologous twists (thereby refining the action of $H_{d-k}$).  Very roughly, introduction of twists allows us to consider cycles in $M^{d+1}$ relative to the boundary, and therefore to ``cut up" the $d-k$-cycles in $H_{d-k}(M,A)$ into more basic pieces.

%The data of a 1-form symmetry action of some abelian group $A$ on some ribbon braided tensor category $\cA$ allows us to simultaneously \emph{twist} $\Sk_A(M)$ by elements of $H_1(M,A)$ and \emph{grade} $\cA$ by elements of $H_1(M,A^\vee)$.  We recall the construction here: 

In order to simplify notation, we will fix $d=3$ and $k=1$ in what follows; however it is straightforward to modify the general notions for arbitrary $d$ and $k$, and also to introduce higher groupoids facilitating twists in higher codimension.

Recall that a simplicial 1-cycle $\gamma$ on $M$ with values in $A$ is by definition an oriented graph embedded into $M$, with edges labelled by elements of $A$, and with the property that at each vertex the sum of all labels of incoming edges equals to the sum of all labels of outgoing edges.

\begin{definition}
Fix an oriented $3$-manifold $M$.  The Picard groupoid $H_{1,2}(M,A)$ has as its objects the abelian group of simplicial $1$-cycles $C_1(M,A)$, with symmetric monoidal product given by addition.  Morphisms from $\gamma$ to $\gamma'$ are homology classes of $2$-chains $\eta$ satisfying $d\eta = \gamma'-\gamma$.
\end{definition}

It follows immediately from the definition that the group of isomorphism classes of objects of $H_{1,2}(M)$ is canonically isomorphic to $H_1(M,A)$, and that the automorphism group of any object is canonically isomorphic to $H_2(M,A)\cong H_1(M,A^\vee)$.  Hence, the data of a functor $H_{1,2}(M,A) \to \Vect$ consists of a vector space $F([\gamma])$ for each $[\gamma]\in H_1(M,A)$, each equipped with an action of he abelian group $H_1(M,A^\vee)$.

\subsection{One-form actions on braided tensor categories.}
Recall from \cite{BJS} that braided tensor categories form into a 4-category $\BrTens$.  This is the 4-category whose:
\begin{itemize}
    \item objects are cp-rigid braided tensor categories,
    \item 1-morphisms from $\cA$ to $\cB$ are $\cA$-$\cB$-central algebras, i.e. cp-rigid tensor categories $\cC$ equipped with braided tensor functors $\cA\otimes \cB^{op}\to \mathcal{Z}(\cC)$,
    \item 2-morphisms from $\cC$ to $\cD$ are $\cC$-$\cD$ bimodule categories, compatibly with the central action,
    \item 3-morphisms are bimodule functors, and
    \item 4-morphisms are natural transformations.
\end{itemize}

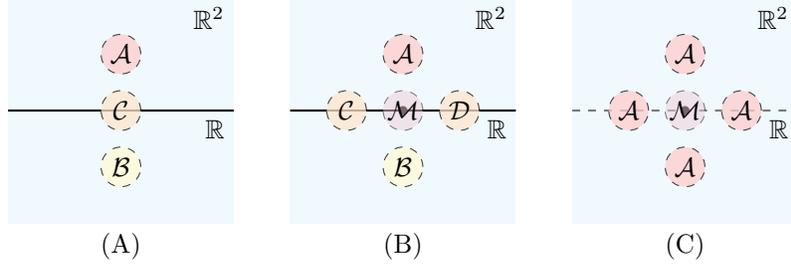
\begin{figure}[!ht]\begin{center}
\begin{tikzpicture}[scale=0.75]

\begin{scope}[shift={(-5,0)}]%basic open set labels
\fill[fill=cyan!5] (-2,-2) rectangle (2,2);
\draw[thick] (-2,0) -- (2,0);
\filldraw[fill=orange!20, dashed,opacity=0.75] (0,0) circle (10pt);
\draw (0,0) node  {$\mathcal{C}$};
\filldraw[fill=red!20, dashed,opacity=0.75] (0,1) circle (10pt);
\draw (0,1) node  {$\cA$};
\filldraw[fill=yellow!20, dashed,opacity=0.75] (0,-1) circle (10pt);
\draw (0,-1) node {$\mathcal{B}$};
\draw (2,2) node [below left] {$\RR^2$};
\draw (2,0) node[below left] {$\RR$};
\draw (0,-2) node[below] {(A)};
\end{scope}

\begin{scope}[shift={(0,0)}]%basic open set labels
\fill[fill=cyan!5] (-2,-2) rectangle (2,2);
\draw[thick] (-2,0) -- (2,0);
\filldraw[fill=black] (0,0) circle (2pt);
\filldraw[fill=orange!20, dashed,opacity=0.75] (-1,0) circle (10pt);
\draw (-1,0) node  {$\mathcal{C}$};
\filldraw[fill=orange!20, dashed,opacity=0.75] (1,0) circle (10pt);
\draw (1,0) node  {$\mathcal{D}$};
\filldraw[fill=purple!20, dashed,opacity=0.5] (0,0) circle (10pt);
\draw (0,0) node  {$\mathcal{M}$};
\filldraw[fill=red!20, dashed,opacity=0.75] (0,1) circle (10pt);
\draw (0,1) node  {$\cA$};
\filldraw[fill=yellow!20, dashed,opacity=0.75] (0,-1) circle (10pt);
\draw (0,-1) node {$\mathcal{B}$};
\draw (2,2) node [below left] {$\RR^2$};
\draw (2,0) node[below left] {$\RR$};
\draw (0,-2) node[below] {(B)};
\end{scope}

\begin{scope}[shift={(5,0)}]%basic open set labels
\fill[fill=cyan!5] (-2,-2) rectangle (2,2);
\draw[thick, dashed, opacity=.5] (-2,0) -- (2,0);
\filldraw[fill=black] (0,0) circle (2pt);
\filldraw[fill=red!20, dashed,opacity=0.75] (-1,0) circle (10pt);
\draw (-1,0) node  {$\mathcal{A}$};
\filldraw[fill=red!20, dashed,opacity=0.75] (1,0) circle (10pt);
\draw (1,0) node  {$\mathcal{A}$};
\filldraw[fill=purple!20, dashed,opacity=0.5] (0,0) circle (10pt);
\draw (0,0) node  {$\mathcal{M}$};
\filldraw[fill=red!20, dashed,opacity=0.75] (0,1) circle (10pt);
\draw (0,1) node  {$\cA$};
\filldraw[fill=red!20, dashed,opacity=0.75] (0,-1) circle (10pt);
\draw (0,-1) node {$\mathcal{A}$};
\draw (2,2) node [below left] {$\RR^2$};
\draw (2,0) node[below left] {$\RR$};
\draw (0,-2) node[below] {(C)};
\end{scope}

\end{tikzpicture}\end{center}
\caption{(A) depicts a pair of braided tensor categories $\cA$ and $\mathcal{B}$, encoded by embeddings of disks $\mathbb{R}^2$, and a 1-morphism between them: an $\cA$-$\mathcal{B}$-central algebra $\mathcal{C}$  encoded by embeddings of intervals $\mathbb{R}$ sandwiched between the $\cA$- and $\mathcal{B}$-regions.  (B) depicts a 2-morphism: a $\mathcal{C}$-$\mathcal{D}$-bimodule structure on $\mathcal{M}$ is encoded by bringing disks into the central point.  (C) depicts the special case $\cA=\mathcal{B}=\mathcal{C}=\mathcal{D}$, in which case we can erase the $\mathbb{R}$, and $\cM$ sits on an isolated point in $\mathbb{R}^2$, and the braided module category structure comes from approaching $\mathcal{M}$ in any direction with $\cA$-disks. }\label{fig:brmodcat}
\end{figure}

In particular, given a braided tensor category $\cA$, its endomorphism 3-category consists of $\cA$-$\cA$-central algebras.  The unit endomorphism is $\cA$ itself with its regular action.  We may then consider the endomorphisms 2-category of $\cA$ thusly regarded, it is precisely the 2-category $\BrMod(\cA)$ of braided module categories we have met previously.  See Figure \ref{fig:brmodcat} for a comparison with the notion in Definition \ref{def:brmodcatdef}

\begin{definition}
A 1-form symmetry action of the abelian group $A$ on the braided tensor category $\cA$ is the data of an $E_2$-monoidal functor $A\to \BrMod(\cA)$, where $A$ is regarded as a discrete symmetric monoidal 2-category.
\end{definition}

More concretely, a 1-form symmetry action of $A$ on $\cA$ is the data of a braided module category $\cM_a\in\BrMod(\cA)$ for each $a\in A$, and fixed equivalences $\phi_{a,b}:M_a\otimes_\cA M_b\to M_{a+b}$ of braided module categories, satisfying an evident pentagon equation.

Here are the primary two sources of examples we will use:

\begin{example}\label{ex:twist-electric}
Let $A$ be any subgroup of the group $Z_{gp}(\cA)$ of group-like Bernstein-central elements of $A$,
    \[
    Z_{gp}(\cA) = \{ a \in \Aut(\id_\cA) \,\, | \,\, a_{V\otimes W}= a_V \otimes a_W\},
    \]
    and define $\cM_a=\cM$ as a category, with $\alpha_{a,m} = (\id_Y\ot a) \circ \sigma_{m,Y}\sigma_{Y,m}$.  A particular case of this example arises when we take $\cA$ to be $\Repq(G_{sc})$, and $A=Z(G)$, the universal grading group of $\Repq(G_{sc})$.  In physical literature (motivated by the Maxwell theory), the resulting action is called the \defterm{electric 1-form symmetry}.
\end{example}
\begin{example}\label{ex:magnetic}
 Let $\widetilde{A}$ be any $A$-graded braided tensor extension of $\cA$, and define $\cM_a$ to be the degree $a$ component $\Repq(\widetilde{G})_a$, with $\alpha$ and $\phi$ the restrictions of the double braiding and the multiplication functor on $\widetilde{A}$.  A particular case of this example arises when we take $\cA$ to be $\Repq(G_{ad})$, $A=\pi_1(G)$, and $\widetilde{A}=\Repq(G_{sc})$.  Again in physical literature, the resulting action is called the \defterm{magnetic 1-form symmetry}.
\end{example}

\subsection{The twisted skein module}
Let us now translate the abstract notions above into the formalism of skein modules, and ultimately linear algebra.

\begin{definition}
Fix a 1-form action $\pi:\cA\to \BrMod(\cA)$.  For each  simplicial cycle $\gamma \in C_1(M,A)$, the \defterm{twisted skein module} $\TwSk(M)[\gamma]$ is the span of all $(\cA,\pi,\gamma)$-labelled ribbon graphs, modulo the $(\cA,\pi,\gamma)$-skein relations.
\end{definition}

Here, a $(\cA,\pi,\gamma)$-labelled ribbon graph is a ribbon graph in $M$, with strands and coupons labelled by objects and morphisms of $\cA$ on the complement of the support of $\gamma$, and with strands and coupons labelled by objects and morphisms of $\cM_a$ along any 1-simplex of $\gamma$ weighted by $a$.  The skein relations are defined in the evident way, by treating separately cylinders which do not intersect the support of $\gamma$, cylinders which intersect a 1-simplex of $\gamma$, and cylinders which intersect a 0-simplex of $\gamma$.  We omit complete definitions.

%Given this data, and given a 1-form symmetry action of $A$, we define $\TwRib(M)[\gamma]$ to be the vector space formally spanned by ribbon graphs which are $\cA$-labelled away from the support of $\cM$, and $\cM_a$-labelled along any 1-simplex labelled $a$.  We then impose the ordinary skein relations on any ball not intersecting the 1-simplex, and we impose $\cM_a$-skein relations on any ball intersecting a single $a$-simplex along a single component.

\begin{proposition}[\cite{GJS2}]
We have a well-defined functor,
\[
\TwSk(M):H_{1,2}(M,Z(G)^\vee)\to \Vect,
\]
mapping each $\gamma\in Z_1(M,Z(G)^\vee)$ to the $\gamma$-twisted skein module, and mapping each $\eta:\gamma_1\to\gamma_2$ to an isomorphism induced by Reshetikhin-Turaev evaluation on a bounding disk.
\end{proposition}

This means simply that for each choice of $\gamma\in C_1(M,A)$ we will define a twisted skein module $\TwSk(M)[\gamma]$, and for each $\eta \in C_2(M,A)/dC_3(M,A)$ with $d\eta = \gamma'-\gamma$ we define an isomorphism,
\[
\TwSk(M)[\eta]: \TwSk(M)[\gamma]\to \TwSk(M)[\gamma'],
\]
and that these isomorphisms will compose compatibly.

Finally, we can consider the two flavours of twisted skein modules coming from Example \ref{ex:SL2PGL2}.

\begin{figure}
    \centering
    \includegraphics[height=6cm]{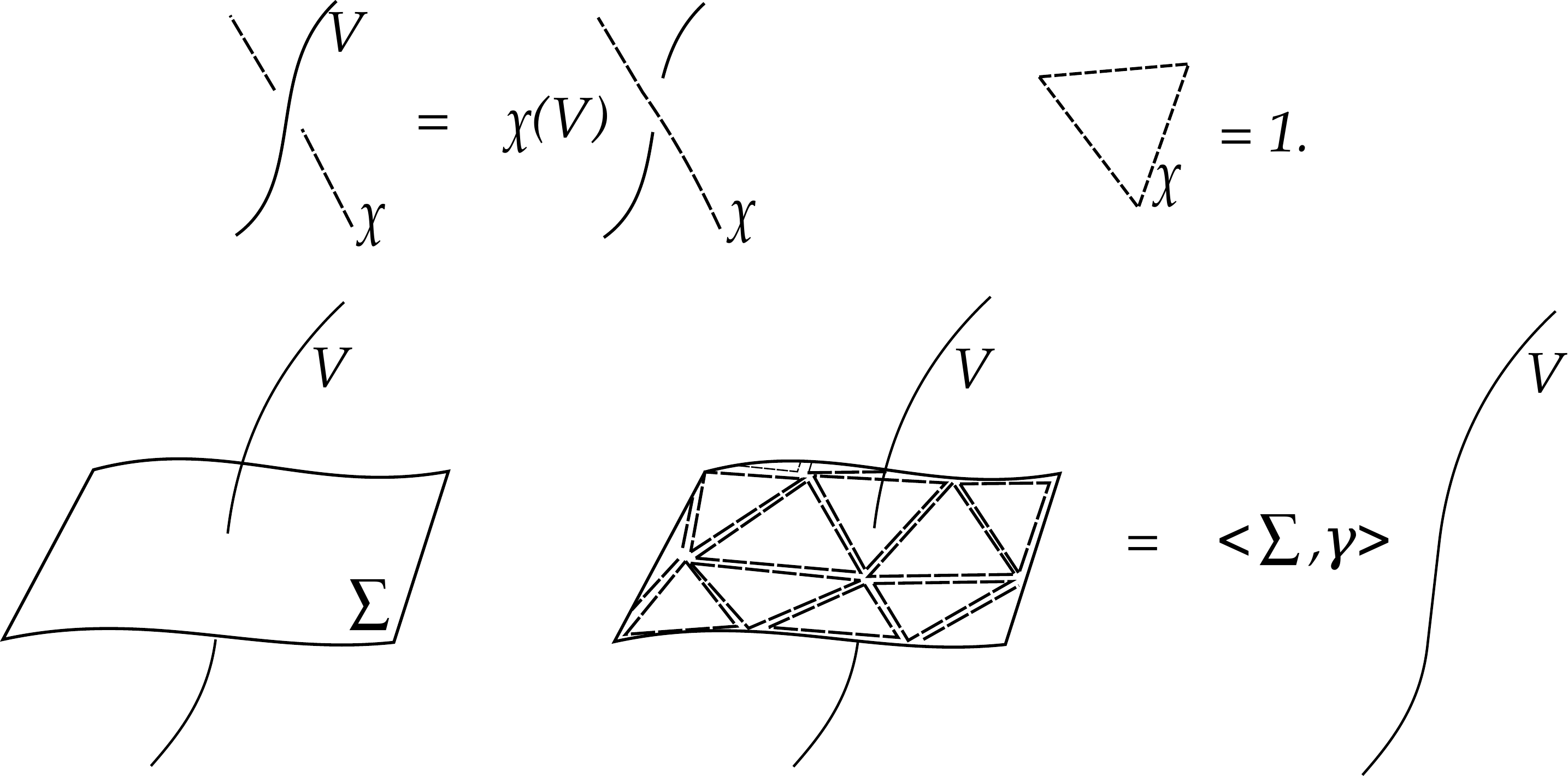}
    \caption{Top left: each character $\chi$ determines a braided module category $\mathcal{M}_\chi$, whose underlying category is $\cA$, and where the double braiding is given by the central action of $\chi$ on each object $V\in\cA$ marking a skein.  Top right: this braided module category is invertible; closed loops contract trivially.  At bottom:  the pairing between some 2-cycle $\Sigma$ and some skein is obtained by first condensing within $\Sigma$ some network of triangles, and contracting each triangle through the skein.  The result is multiplication by the Poincare pairing $\langle \Sigma,\gamma\rangle$ between the cycle $\Sigma$ and the support 1-cycle $\gamma$ of the skein.  We note that invertibility of $\chi$ implies already that adjacent edges of the triangulation cancel pairwise, hence the resulting skein module is not twisted despite the introduction of defects, so that we indeed obtain an endomorphism.}
    \label{fig:H2fig}
\end{figure}

\begin{example} Let $A=Z(G)$.  We have seen in Example \ref{ex:twist-electric} that we have an action of $Z(G)$ on $\Repq(G)$ by 1-form symmetries; the resulting action by $H_2(M,Z(G))$ on $\Sk(M) = \TwSk(M)[0]$ allows us to decompose $\Sk(M)$ according to the characters $H_2(M,Z(G))^\vee \cong H_1(M,Z(G)^\vee)$.  The skein-theoretic interpretation is depicted in Figure \ref{fig:H2fig}.  The resulting graded decomposition is the elementary one observed at the start of Section \ref{sec:gradingstwistings}.  This grading, in its elementary form, has featured prominently in the computations of \cite{Carrega}, \cite{Gilmer}, \cite{GilmerMasbaum} and \cite{Detcherry-Wolff}.

To our knowledge, the twisted skein modules $\TwSk(M)[\gamma]$ for non-trivial $\gamma\in H_1(M,Z(G))$ are new, and so do not have an {\it a priori} elementary description.  However, it is easy to describe the twists in this case:  the condition on defects says that when a 1-simplex of a ribbon graph labelled with a simple representation $V$ passes through a 1-simplex of the defect with weight $z\in Z(G)$, it picks up a scalar factor equal to the evaluation of $z$ on $V$.\end{example}

\begin{example}Let $A=\pi_1(G)$.  We have seen in Example $\ref{ex:magnetic}$ that we obtain an action of $\pi_1(G)$ on $\Repq(G)$ by 1-form symmetries.  To our knowledge both the resulting action of $H_2(M,Z(G))$ and hence the resulting decomposition of $H_1(M,\pi_1(G)^\vee)$ on $\Sk(M)$, and the twisted skein modules are new, and so do not have an {\it a priori} elementary description.  However, it is again easy to define twists in this case.  Given a $G^{sc}$-labelled ribbon graph with simple labels, we can produce a 1-simplex in $C_1(M,\pi_1(G))$, by replacing each simple label $V$ by its grading with respect to $\pi_1(G)$, regarding $\Repq(G^{sc})$ as a $\pi_1(G)$-graded extension of $\Repq(G)$.  The space $\TwSk(M)[\gamma]$ is the span of all $G^{sc}$-labelled ribbon graphs, whose underlying $1$-cocycle equals $\gamma$, and modulo only those skein relations also with underlying $1$-cocycle $\gamma$.
\end{example}

For any group $G$, and for $A=Z$ or $\pi_1$, and for $a\in H_1(M,A^\vee)$ and $b\in H_1(M,A)$ we denote by $\Sk_{G,A}^{a,b}(M)$ the $a$-graded component of the $b$-twisted skein module.  It follows from the definition that we have a natural linear map of vector spaces,
\[
\Phi: \Sk_{G^{ad},Z}(M) \to \Sk_{G^{sc}}(M),
\]
which regards a $\Repq(G^{ad})$-labelled ribbon graph as a $\Repq(G^{sc})$-labelled ribbon graph.  An important observation is that while $\Repq(G^{ad})\to \Repq(G^{sc})$ is fully faithful (so that in a ball this map is simply an isomorphism), it is \emph{not} an isomorphism globally, and has both kernel and cokernel.  The cokernel is clear: by construction the image of $\Phi$ is contained in the degree zero part of $\Sk_{G^{sc}}(M).$  The kernel is more subtle, but it stems from the fact that on the source, by definition, only skein relations which have trivial underlying 1-cycle are imposed, whereas in the target they are not.  In the skein module $\Sk_{G^{sc}}(M)$, two ribbon graphs of degree zero can be related in a way that nevertheless involves elements of nontrivial degree.  This is captured cleanly in the following

\begin{theorem}[\cite{GJS2}]\label{thm:easy-identity} Let $M$ be an oriented $3$-manifold (not necessarily closed), and let $q\in\mathbb{C}^\times$ be an arbitrary complex number.  Then we have a canonical isomorphism,
\[
\Sk_{G/Z(G),\pi_1}^{a,b}(M) \cong \Sk_{G,Z}^{a,b}(M).
\]
\end{theorem}

Unpacking, it says that the map $\Phi$ above actually extends firstly to a map
\[
\Sk_{G^{ad},\pi_1}(M) \to \bigoplus_{b\in H_1(M,Z(G))} \Sk^{0,b}_{G^{ad},Z}(M),
\]
which is injective, and then again to a map,
\[
\bigoplus_{a,b}\Sk^{a,b}_{G^{sc},\pi_1}(M) \to \bigoplus_{a,b} \Sk^{a,b}_{G^{sc},Z}(M),
\]
which is an isomorphism.  We note that the indices $a$ and $b$ appear on botht he LHS and RHS above.  This is justified by the natural isomorphisms, $\pi_1(G^{ad})\cong Z(G^{sc})^\vee$.

The reader familiar with the language of gauging of higher form symmetries may not find this surprising:  Indeed, we have $G^{ad} = G/Z(G)$, hence $\Repq(G^{ad})$ is a gauging of the electric $1$-form symmetry by $Z(G)$; gauging an electric 1-form symmetry induces on the new theory a dual magnetic 1-form symmetry, with an isomorphism like that above.

\section{The bigraded refinement of the skein conjecture}
\label{sec:bigraded-refinement}
In total, we have defined bi-graded skein module,
\[
\Sk_{G^{sc},Z}^{a,b}(M), \qquad \textrm{for $a\in H_1(M,Z(G^{sc})^\vee), b\in H_1(M,Z(G^{sc}))$},
\]
using the electric 1-form symmetry, and
\[
\Sk_{G^{ad},\pi_1}^{a,b}(M), \qquad \textrm{for $a\in H_1(M,\pi_1), b\in H_1(M,\pi_1(G^{ad})^\vee)$},
\]
using the magnetic 1-form symmetry, where in each case $a$ denotes the grading degree, and $b$ denotes the twist.  Now recall that under Langlands/S-duality, the electric and magnetic 1-form symmetries are interchanged.  Thus the same electric-magnetic duality principle which informs Conjecture \ref{conj:main-conj} also predicts:

\begin{conjecture}[Bigraded refinement of skeins conjecture]\label{conj:bigraded-refinement}
Suppose that $q\in\mathbb{C}^\times$ is transcendental and that $M$ is a closed oriented $3$-manifold.  Then we have an equality of dimensions,
\[
\dim\Sk_{G,Z}^{a,b}(M) = \dim\Sk_{{^LG},\pi_1}^{b,a}(M).
\]
\end{conjecture}

Now suppose that $G$ is simply connected of ADE type, so that ${{}^L\mathfrak{g}}=\mathfrak{g}$.  In this case the Langlands dual group ${^LG}=G^{ad}$ has a second manifestation as a quotient $G/Z(G)$ of $G$ by the electric 1-form symmetry we exploited to define the bi-grading.  Hence we may apply Theorem \ref{thm:easy-identity}, to obtain the following equivalent

\begin{conjecture}[Reformulation of refined skein conjecture]\label{conj:reformulated-refined} Suppose that $G$ is simply connected of ADE type, that $q\in \mathbb{C}^\times$ is transcendental and that $M$ is a closed oriented $3$-manifold.  Then we have an equality of dimensions,
\[
\dim\Sk_{G,Z}^{a,b}(M) = \dim\Sk_{G,Z}^{b,a}(M),
\]
for all $a\in H_1(M,Z(G)^{\vee}), b \in H_1(M,Z(G))$.
\end{conjecture}

\begin{remark} The reformulation of Conjecture \ref{conj:bigraded-refinement} is very convenient:  applying it to the case $G=\SL_N$, it means we can check the refined conjecture completely internally to either the $Z$-twists on $\SL_N$ or alternatively the $\pi_1$-twists on $\PGL_N$.  Because more is known about $\SL_N$ skein theory, we work with the $Z$-twists in the next section.
\end{remark}

\begin{remark}
The isomorphism in Theorem \ref{thm:easy-identity} is completely canonical, holds for arbitrary $q$ and also for 3-manifolds with boundary: it is a very general principle, and in particular, there was an a priori natural alignment of the indices $a$ and $b$ appearing on either side of the asserted isomorphism.  By contrast, the indices $a$ and $b$ do not themselves live in the same group naturally, but rather live in Pontryagin dual groups. Hence Conjecture \ref{conj:reformulated-refined} seems not to parse semantically.

Indeed, to make sense of Conjecture \ref{conj:reformulated-refined}, one must choose an isomorphism $Z(G) \cong Z(G)^\vee$, though clearly the statement itself does not depend on the choice.  This should be regarded as a further symptom of the fact that Langlands duality should be swapping $Z(G)$ with $\pi_1({^LG})^\vee$, but that we are identifying $\pi_1({^LG})^\vee$ with $Z(G)$ in the ADE case.
\end{remark}

Recall that the untwisted skein modules in Conjecture \ref{conj:main-conj} are obtained as sums over the untwisted graded components.  Hence, Conjecture \ref{conj:main-conj} follows from the following weakening of Conjecture \ref{conj:reformulated-refined}.

\begin{conjecture}\label{conj:special-case} 
Suppose that $G$ is simply connected of ADE type, that $q\in \mathbb{C}^\times$ is transcendental and that $M$ is a closed oriented $3$-manifold.  Then we have an equality of dimensions,
\[
\dim \Sk_{G,Z}^{a,0}(M) = \dim \Sk_{G,Z}^{0,a}(M),
\]
for all $a\in H_1(M,Z(G)^{\vee}$.
\end{conjecture}

\section{Some evidence for the conjecture}
\label{sec:evidence}
Dimensions of skein modules are very difficult to compute in elementary terms.  Indeed, this is underscored by the twenty five years which passed between the two fundamental examples $S^2\times S^1$ \cite{HostePrzytycki} and $T^3$ \cite{Carrega,Gilmer}.  Since then a dam has broken and these dimensions are being computed using a combination of elementary and categorical techniques.  This allows us now to prove our conjecture in a number of special cases.
%
%\begin{remark}
%To avoid confusion, we note that there is no material difference between the skein theory defined by a semi-simple affine algebraic group $G$ over $\mathbb{C}$ and its maximal compact Lie subgroup $K$, as these have the same finite-dimensional representation theory.  For instance, it is equivalent to talk about $SL_2$- and $SU(2)$-skeins, and likewise to talk about $PGL_2$ and $SO(3)$-skeins.
%\end{remark}
%
Recall that the Lie algebra $\mathfrak{sl}_N$ is fixed under Langlands duality; however on the group level, Langlands duality interchanges the simply connected form with the adjoint form on the dual group, so it is meaningful to contemplate Langlands duality for the pair ($\SL_N$, $\PGL_N$).

We focus on two families of skein modules, where complete computations can be done.  In the first family, we fix the Langlands dual pair to be ($\SL_2$, $\PGL_2$), and consider mapping cylinders $\Sigma_g\times S^1$ of arbitrary genus $g$, while in the second family we instead fix the 3-manifold to be $T^3$, and vary the Langlands dual pair to be ($\SL_N$, $\PGL_N$). For the first family, we recall the following:

\begin{theorem}[\cite{GilmerMasbaum},  \cite{Detcherry-Wolff}]
Suppose that the quantum parameter is generic.  Then we have:
\[
\dim \Sk_{SL_2}(\Sigma_g\times S^1) = 2^{2g+1} + 2g-1.
\]
With respect to the electric 1-form symmetry action, we have:
\[
\dim \Sk_{\SL_2}^{a,0}(\Sigma_g\times S^1) = \left\{\begin{array}{ll} g+1, &a=(0,0)\\ g, & a=(1,0)\\ 1, & \textrm{otherwise}\end{array}\right.,
\]
where the cases refer to the K\"unneth decomposition,
\begin{align*}
    H_1(\Sigma_g\times S^1,\mathbb{Z}/2) &= H_0(\Sigma,\mathbb{Z}/2)\otimes H_1(S^1) \oplus H_1(\Sigma,\mathbb{Z}/2)\otimes H_0(S^1)\\
    &= H_0(\Sigma,\mathbb{Z}/2)\oplus H_1(\Sigma,\mathbb{Z}/2).
\end{align*}
\end{theorem}

The first non-trivial confirmation of Langlands duality for skein module is then given by the following

\begin{theorem}[\cite{GJS2}, to appear]\label{thm:Sigmag}
Suppose that the quantum parameter is generic.  Then we have:
\[
\dim \Sk_{\PGL_2}(\Sigma_g\times S^1) = 2^{2g+1} + 2g-1
\]
\end{theorem}

In fact the proof is not by computing with $\Sk_{\PGL_2}$, but instead by applying Corollary \ref{conj:special-case}.  More precisely we prove:
\[
\dim \Sk_{\SL_2}^{0,a}(\Sigma_g\times S^1) = \left\{\begin{array}{ll} g+1, &a=(0,0)\\ g, & a=(1,0)\\ 1, & \textrm{otherwise}\end{array}\right..
\]
The first case is, of course, already proved by Carrega and Gilmer, we repeat it only for emphasis.  For the second case, we use the fact that the twisted skein category at $q=1$ becomes the category of quasi-coherent sheaves on the \emph{twisted character variety}.  This is in fact a \emph{smooth affine variety}, hence the Hochschild homology of its quantization to the skein category is identified with the middle homology group of the twisted character variety, which was shown to be $g$-dimensional by Hausel--Rodriguez-Villegas \cite{HRV}, following Hitchin.

Turning now to the second family of skein modules, we recall the following notation: $\mathcal{P}(d)$ denotes the partition number of $d$, i.e. the number of ways to write $d$ as a sum of non-negative integers, $J_3=C \ast \mu$ is the third Jordan totient function, obtained as the M\"obius inversion of the cube function $C(n)=n^3$.

\begin{theorem}[\cite{GJVY}, to appear]
Suppose that the quantum parameter is generic.  Then we have:
\[
\dim \Sk_{SL_N}(T^3) = (\mathcal{P}\ast J_3)(N) = \sum_{N=d\cdot e\cdot f}\mathcal{P}(d)\cdot e^3 \cdot \mu(f).
\]
With respect to the electric 1-form symmetry action we have
\[
\Sk_{SL_N}(T^3)^{(a,b,c),0} = \operatorname{gcd(a,b,c,N)}, 
\]
where we write $(a,b,c)$ for an element of $(\mathbb{Z}/N)^3 = H_1(T^3,\mathbb{Z}/N)$.
\end{theorem}

This theorem is proved by first obtaining an equivalence between the skein category of the torus and the category of modules for a double affine Hecke algebra, subsequently giving an isomorphism between this double affine Hecke algebra and a smash product of a quantum torus and the symmetric group, and finally computing the Hochschild homology directly by combinatorial means.

Our second non-trivial confirmation of Langlands duality is given by the following

\begin{theorem}[\cite{GJS2}]
Suppose that the quantum parameter is generic, and that $N$ is prime.  Then we have:
\[
\dim \Sk_{PGL_N}(T^3) = (\mathcal{P}\ast J_3)(N) = \sum_{N=d\cdot e\cdot f}\mathcal{P}(d)\cdot e^3 \cdot \mu(f).
\]
With respect to the magnetic 1-form symmetry action we have
\[
\Sk_{\PGL_N}(T^3)^{(a,b,c),0} = \operatorname{gcd(a,b,c,N)}=1. 
\]
\end{theorem}

Our proof is similar to that of Theorem \ref{thm:Sigmag}.  We apply Conjecture \ref{conj:special-case}, and instead prove:
\[
\dim \Sk_{\SL_N}^{0,(a,b,c)}(T^3) = 1,\quad \textrm{if $(a,b,c)\neq 0$}.
\]
For this, we rely again on the fact that the twisted character variety is smooth, and that its middle homology is known to be one-dimensional.  Of course, we expect the theorem to hold when $N$ is not prime, but in this case the twisted character variety is not smooth, so we can't easily apply the same techniques for the proof.

\printbibliography

\end{document}

%% file: macros.tex
\newcommand{\cA}{\mathcal{A}}

\newcommand{\C}{\mathbf{C}}

\newcommand{\D}{\mathcal{D}}
\newcommand{\Sk}{\mathrm{Sk}}

\newcommand{\Coh}{\mathrm{Coh}}
\newcommand{\SL}{\mathrm{SL}}
\newcommand{\PGL}{\mathrm{PGL}}
\newcommand{\DD}{\mathbb{D}^2}
\newcommand{\SkCat}{\operatorname{SkCat}}
\newcommand{\TwSk}{\operatorname{TwSk}}
\newcommand{\Vect}{\mathrm{Vect}}
\newcommand{\BrTens}{\mathrm{BrTens}}
\newcommand{\rev}{\operatorname{rev}}

\newcommand{\Bun}{\mathrm{Bun}}
\newcommand{\CP}{\mathbf{CP}}

\newcommand{\Hom}{\mathrm{Hom}}

\newcommand{\Loc}{\mathrm{Loc}}
\newcommand{\Ch}{\mathrm{Ch}}

\newcommand{\Rep}{\mathrm{Rep}}
\newcommand{\Repq}{\mathrm{Rep_q}}
\newcommand{\Shv}{\mathrm{Shv}}
\newcommand{\BrMod}{\mathrm{BrMod}}
\newcommand{\Aut}{\mathrm{Aut}}
\newcommand{\id}{\mathrm{id}}
\newcommand{\ot}{\otimes}
\newcommand{\cD}{\mathcal{D}}
\newcommand{\cC}{\mathcal{C}}
\newcommand{\cB}{\mathcal{B}}
\newcommand{\cM}{\mathcal{M}}
\newcommand{\RR}{\mathbb{R}}

\numberwithin{equation}{section}
\newtheorem*{theorem*}{Theorem}
\newtheorem{theorem}{Theorem}[section]
\newtheorem{theorem/def}[theorem]{Theorem/Definition}

\newtheorem{proposition}[theorem]{Proposition}
\newtheorem{conjecture}[theorem]{Conjecture}

\theoremstyle{definition}
\newtheorem{definition}[theorem]{Definition}
\newtheorem{remark}[theorem]{Remark}

\newtheorem{example}[theorem]{Example}
\newtheorem{def/prop}[theorem]{Definition/Proposition}
\newtheorem{question}[theorem]{Question}

\newcommand{\defterm}[1]{\textbf{#1}}